\documentclass[a4paper]{article}

\newif
\ifdebug
\debugtrue 

\usepackage[unicode]{hyperref}
\usepackage{amsfonts}
\usepackage{amssymb}
\usepackage{amsmath}
\usepackage{euscript}
\usepackage{amsthm}   
\usepackage{amsfonts} 
\usepackage[matrix,arrow,curve]{xy}\CompileMatrices

\ifdebug 
    \usepackage{xcolor} 
    \def\comment#1{\textcolor{gray!80}{#1}} 
    \def\fcomment#1{\footnote{\textcolor{gray!80}{#1}}}
\else 
    \def\comment#1{} 
    \def\fcomment#1{}
\fi

\newcommand{\vxymatrix}[1]{\vcenter{\xymatrix{#1}}}

\newcommand{\obbR}{{\mbox{\kern1.4pt \raisebox{7pt} {\scriptsize$\circ$} \kern-9.1pt}\mathbb{R}}}
\newcommand{\indobbR}{{\mbox{\kern0.7pt \raisebox{5pt} {\tiny$\circ$} \kern-7.8pt}\mathbb{R}}}
\newcommand{\mbbR}{{\mathbb{R}}}

\newcommand{\be}{\begin{equation}}
\newcommand{\ee}[1]{\label{#1}\end{equation}}

\theoremstyle{plain}
\theoremstyle{plain}
\newtheorem{thm}{Theorem}[section]
\newtheorem{lm}[thm]{Lemma}

\theoremstyle{definition}
\newtheorem{definition}[thm]{Definition}
\theoremstyle{remark}

\DeclareMathOperator{\Cl}{Cl}
\DeclareMathOperator{\Int}{Int}
\DeclareMathOperator{\Diff}{Diff}

\DeclareMathOperator{\Spec}{Spec}

\DeclareMathOperator{\Smbl}{Smbl}

\begin{document}
"Poisson structure on manifolds with singularities"
\\
M. Sorokina

\section{\Large Introduction}
The configuration spaces of many real mechanical systems appear to be manifolds with singularity. A singularity often indicates that geometry of motion might change at the point.

In applied problems of mechanics, these singularities usually cause heavy loads on system parts and, due to inevitable miscalculations, make the system uncontrollable, so they should be reduced while solving. On the other hand, we face the conceptual difficulties even describing mechanics of the ideal models. Namely, since the configuration space is not a smooth manifold the whole set of techniques of Hamiltonian mechanics cannot be applied.

In this paper we present a way of conquering the aforementioned conceptual problem by considering a certain algebra whose real spectrum coincides with the above configuration space. The structure of this algebra is completely determined by the geometry of the singularity. For a broad class of singularities, the desired algebra can be described directly as the pullback of the two already known algebras. Availability of the algebra enables us to use the Differential operator theory.

A role of the phase space, i. e. a space on which Hamiltonian Formalism is being given, is being played by differential operator symbol algebra which is explicitly declared in the sections.

Flat linkages are simple examples of mechanical systems to which this algorithm is applicable.

Below we build a Poisson structure on a manifold with a one-dimensional singularity. Similar result can be obtained for some other kinds of singularities.

As an example, there are two intersecting manifolds. In this case, algebra of functions at the intersection is being given depending on how the manifolds had been initially given. We can consider that the manifolds lie in a many-dimensional space and unite there. Then, algebra of functions at the union is being obtained by restriction of algebra of all functions on this manifold space to this union. Algebra at the intersection is the algebra of equivalence classes of functions at the sewing: two functions are equivalent if their restrictions to the first component match, if restrictions to the second component match or they can be connected with a finite number of those equivalences.

In particular, at the end of this report we present the results for the specific case of a configuration space consisting of two curves on a plane having arbitrary
order of contact.

The author is grateful to A. M. Vinogradov who suggested the problem considered and for plenty of helpful discussions.

\newpage
\subsection{\Large Definitions}
Let $A$ will always be a commutative, associative algebra with unit over $\mbbR$.
All $\mbbR$-algebra homomorphisms $h\colon A_1\to A_2$ are assumed unital (i.e., $h$ maps the identity element of $A_1$ to the identity element of $A_2$).
\begin{definition}
Let $A$ be a $\mbbR$-algebra. Then a $\mbbR$-homomorphism $\Delta\colon A\to A$ is called a \textit{linear differential operator of order $\leq k$} with values in $A$ if for any $a_0,...a_k \in A$ we have the identity
$$(\delta_{a_0}\circ\dots\circ\delta_{a_k})(\Delta)=0.$$
\end{definition}
Let us denote the set of all differential operator of order $\leq k$ acting from $A$ to $A$ by $\Diff_{k}$. This set is stable with respect to summation and multiplication by elements of the algebra $A$. Therefore, it is naturally endowed with an $A$-module structure.

\begin{definition}
Embedding of $A$-modules $\Diff_{k-1}(A)\subset\Diff_k(A)$ allows one to define the quotient module $S_k(A)=\frac{\Diff_k(A)}{\Diff_{k-1}(A)},$ which is called \textit{the module of symbols of order $k$ (or the module of $k$-symbols)}. 
The coset of an operator $\Delta\in\Diff_k(A)$ modulo $\Diff_{k-1}(A)$ will be
denoted by $smbl_k\Delta$ and called \textit{the symbol of $\Delta$.}
Let us define the \textit{algebra of symbols} for the algebra $A$ by setting
$$Smbl(A)= \bigoplus_{n=0}^{\infty} S_n(A).$$
\end{definition}
The operation of multiplication in $\Smbl(A)$ is induced by the composition of
differential operators. To be more precise, for two elements
$$smbl_l(\Delta_1)\in S_l(A), smbl_m(\Delta_2)\in S_m(A)$$
let us set by definition 
$$smbl_l(\Delta_1)\cdot smbl_m(\Delta_2):=smbl_{l+m}(\Delta_1\circ\Delta_2)\in S_{l+m}(A).$$
This operation is well defined, since the result does not depend on the choice of representatives in the cosets $smbl_l(\Delta_1)$ and $smbl_m(\Delta_2).$
$\Smbl(A)$ is a commutative graded $\mbbR$-algebra.
\begin{definition}
Now let $\Delta\in\Diff_m(A)$ and $\nabla\in\Diff_l(A).$ Then $[\Delta,\nabla]\in\Diff_{m+l-1}(A).$
One can assign to the pair $(smbl_m(\Delta),smbl_l(\nabla))$ the element
$$\{smbl_m(\Delta),smbl_l(\nabla)\} = smbl_{l+m-1}([\Delta,\nabla])\in S_{l+m-1}(A),$$
which is well defined, i.e., does not depend on the choice of representatives
in the cosets.
Denote by $\{.,.\}$ any \textit{Poisson bracket.} operation. The operation $\{.,.\}$ is $\mbbR$-linear and skew-symmetric. It satisfies the Jacobi identity, since the commutator of linear differential operators satisfies this identity. Thus, $Smbl(A)$ is a Lie algebra with respect to this operation.

\end{definition}
Any manifold $M$ is determined by the smooth $\mbbR$-algebra $A=C^{\infty}(M)$ of functions on it, each point $x$ on $M$ being the $\mbbR$-algebra homomorphism $\colon A\to \mbbR$ that assigns to every function $f \in A$ its value $f(x)$ at the point $x$.

Let $f\colon M\to N$ be a map between manifolds then the map $$f^*\colon C^\infty(N)\to C^\infty\\\\(M)\ \colon a\mapsto a\circ f$$ is a homomorphism of corresponding algebras. 
Let $A$ be a commutative algebra with unit. We define  $|A|=\Spec_{\mbbR} A$ as a set of all unitary homomorphisms from $A$ to $\mbbR$.

Thus, $\Spec_{\mbbR} A$  is a spectrum of algebra $A$, homeomorphism $h$ is the point of algebra.
Let us equip $\Spec_{\mbbR} A$ with Zarissky topology generated by the sets $U_a=\{h\in A\mid h(a)\neq 0\}$.

For a manifold $M$ there is a natural homeomorphism
$$\varphi\colon M \to \Spec(C^\infty(M))\ \colon x \mapsto h_x,$$ where $h_x\in\Spec(C^\infty(M))$ is given by $h_x(a)=a(x)$.

\begin{definition}
Let us fix a coordinate system $x_1,\dots, x_n$ in a neighborhood $U$ of a point
$z$. Recall that a domain U is called \textit{starlike} with respect to $z$ if together with any point $y\in U$ it contains the whole interval $(z,y)$.
\end{definition}

\begin{definition}
\textit{Hadamard’s lemma:} Any smooth function $f$ in a starlike neighborhood of a point $z$ is representable in the form
$$f(x)=f(z)+ \sum_{i=1}^n (x_i-z_i)g_i(x), g_i\in \ {C^{\infty}(U)},$$
where $g_i$ are smooth functions.
\end{definition}

If the functions $f,g\in C^{\infty}(M)$ concide in some neighbourhood of $U\ni z$ then for any differential operator $\Delta$ we have ${{\Delta (f)}|}_U={{\Delta (g)}|}_U$. Hence for any differential
operator $\Delta\in\Diff_k(C^{\infty}(M))$ the restriction 
$${\Delta|}_U \colon C^{\infty}(U)\to C^{\infty}(U)$$
on any open domain $U \subset M$ is correctly defined. We call this \textit{locality principle.}

\newpage
\section{\Large Pullback of algebras}

Let us consider objects $X, Y$ и $B$ and morphisms $\phi\colon X\to B$ и $\nu\colon Y\to B$

 \textit{Pullback} is an object $Z$ with morphisms $\alpha\colon Z\to X$, $\beta\colon Z\to Y$ such that $\nu\circ\beta = \phi\circ\alpha$ and for any object  $W$ and morphisms $$\alpha'\colon W\to X,\ \beta'\colon W\to Y$$ from the equality  $\nu\circ\beta' = \phi\circ\alpha'$ implies  that there must exist a unique morphism $\gamma\colon W\to Z$ such that $\alpha' = \alpha\circ\gamma$ and $\beta' = \beta\circ\gamma$. The pullback, if it exists, is unique up to a unique isomorphism. $[7].$             
$$
\xymatrix{
  W\ar@/_/[ddr]_{\alpha'} \ar@/^/[drr]^{\beta'}
   \ar@{.>}[dr]|-{\gamma}            \\
  &  Z \ar[d]^{\alpha} \ar[r]_{\beta}
                 & Y \ar[d]_{\nu}       \\
  & X \ar[r]^{\phi}   & B                }
$$
In the category of sets, a pullback of f and g is given by the set:
$$X\times_{B}Y = \{(x,y)\in X\times Y\mid\phi(x)=\nu(y) \} .$$

Let us consider algebras $A_1 = C^\infty(M_1)$, $A_2 = C^\infty(M_2)$ and homomorphisms $p_i\colon A_i\to C$, $C$ - certain $\mbbR$-algebra. 
\\    
Let us consider pullback:  
    $$\xymatrix{A\ar@{->>}[r]^{\pi_2}\ar@{->>}[d]_{\pi_1} & A_2\ar@{->>}[d]^{p_2}\\ A_1 \ar@{->>}[r]_{p_1}&C}$$
$A=A_1 \oplus_C A_2 =\{(a_1,a_2)\in A_1\oplus A_2: p_1(a_1) = p_2(a_2)\}.$    
   
\textbf{Example 1.} 
$A_1=A_2=C^\infty(\mbbR)$, $C=\mbbR$, $p_1(f)=f(0),\,p_2(g)=g(0)$.

\medskip

\textbf{Example 2.} 
$A_1=A_2=C^\infty(\mbbR)$, $C=\mbbR[x]/(x^{k+1}),$
$$p_i(f)=\sum_{n=0}^k \frac{f^{(n)}(0) \epsilon^n}{n!}$$

\newpage

\subsection{\Large Spectrum of the algebra}
 Let us consider the commutative diagram. 
\\    
    $\mbbR$-algebras homomorphisms  $h_{A_i}$  are defined as a composition 
    $h_{A_i}=h_C \circ p_i$. 
$$
\xymatrix{ A\ar@{->>}[r]^{\pi_2}\ar@{->>}[d]_{\pi_1} & A_2\ar@{->>}[d]^{p_2}\\ A_1 \ar@{->>}[r]_{p_1}& C\\ & &R\ar@{<-}[ul]|-{h_C}\ar@{<-}[uul]_{h_{A_2}} \ar@{<-}[ull]^{h_{A_1}} \ar@{<.}[ulul]}
$$
There are two maps of spectrums $i_{1,2}\colon|C|\to|A_i|$ injective since homomorphisms $p_i$ are epimorphisms. It means that $|C|$ may be considered as a subset of $|A_i|,$ and it makes sense to speak about splice of spectrums.
\begin{lm}
Spectrum of algebra $A$ is equal to splice of spectrums of algebras $A_1$ and $A_2$ by spectrum of algebra $C$:
$$\Spec_\mbbR{A}=\frac{\Spec_\mbbR{A_1}\coprod \Spec_\mbbR{A_2}}{\Spec_\mbbR{C}}$$
\end{lm}
\begin{proof}
Let $h\in |A|$. First let us prove that $h$ lies in at least one of the images 
$|\pi_i|(|A_i|)$ ($i=1,2$), i.e. that at least one of the 
$\ker \pi_i$ is contained in  $\ker h$. Assume this is not true.
Note that any element of $\ker\pi_1$ is a pair $(0,a_2)$ and any element of
$\ker\pi_2$ is a pair $(a_1,0)$. Thus, there exist $a_1\in A_1$ and $a_2\in A_2$ such that $h(a_1,0)\neq0$ and $h(0,a_2)\neq 0$. Then $$h(0,0)=h(a_1,0) h(0,a_2)\neq 0,$$ 
which contradicts to the fact that $h$ is a homomoprhism. 

Now let's investgate those $h$ which belong to both images. Assume 
$$h(a_1,a_2) = h_1(a_1) = h_2(a_2).$$ 
It is easy to conclude that  $h$ is uniquely determined by the value  $p_1(a_1)=p_2(a_2)$.
Therefore $h$ lies in $|p_1\circ\pi_1|(|C|)$.  Obviously the converse is true too, i.e.  $h$ lies in both images iff it lies in $|p_1\circ\pi_1|(|C|)$.  
Hence the set $|A|$ is indeed a gluing of $|A_1|$ and $|A_2|$ along $|C|$.

To prove that the topology on  $|A|$ is the gluing topology we need to show that a $U\subset |A|$ is closed iff it's preimages $|\pi_i|^{-1}(U)$ are closed ($i=1,2$). 
This follows from the fact that $|\pi_i|$ are continous and closed (since $\pi_i$ is a surjective homomorphism). 
\end{proof}


\newpage
\subsection{\Large Construction of differential operators}
This section contains explicit description of differential operators of an algebra $A$ which is a cartesian square of the algebras $A_1, A_2$ over the algebra $C$. 
Let us assume that for $i=1,2$ the following condition holds: 
$$\Cl\Int |\pi_i|(|A_i|) = |\pi_i|(|A_i|).$$ 
\begin{lm}
For every differential operator $\Delta\in\Diff_k(A)$ there exists a unique linear map $\Delta_1\colon A_1\to A_1$ which makes the diagram  \eqref{d1} commutative.
This linear map is a differential operator and lies in $\Diff_k(A_1).$ (Analogous statement is true for $\Delta_2$.)
\end{lm}
\be\vxymatrix{
A\ar[r]^\Delta\ar@{->>}[d]^{\pi_1} &A\ar@{->>}[d]^{\pi_1}\\
A_1\ar[r]^{\Delta_1} &A_1
}\ee{d1}
\begin{proof}
Let us show that $\pi_1(\Delta \ker\pi_1)=\{0\}$.

Given an open subset $U$ of the spectrum of a $C^{\infty}$-closed alebra $A$, $\varphi\in A$ and $\Delta\in\Diff^k(A)$
it follows from the locality principle that the equality  ${\phi|}_U=0$ implies ${(\Delta\phi)|}_U=0$.

Given $U=\Int |\pi_1|(|A_1|)\subset |A|$ let us take $\varphi \in \ker \pi_1$. Then ${\varphi|}_U=0$ hence ${(\Delta\varphi)|}_U=0$. Therefore ${(\Delta\varphi)|}_{\Cl U} = 0$ and hence
$\Delta\varphi\in\ker \pi_1$. Thus we proved that $\pi_1(\Delta \ker\pi_1)=\{0\}$.

It follows from above that the linear map $\Delta_1$ such that the diagram \eqref{d1} is commutative exists and is uniquely defined.

Now let us show that it is a differential operator of the order less or equal to $k$. Let us take  $k+1$ element $a_j\in A_1$
($j=0,...,k$). For each of these let us choose $\tilde a_j
\in \pi_1^{-1}(a_j)$. Then as we know
$$\left((\delta_{\tilde a_0}\circ\dots\circ\delta_{\tilde a_k})(\Delta)\right) \varphi=0.$$
Applying $\pi_1$ to the above equality and using
\begin{align}
\pi_1\circ\Delta = \Delta_1 \circ \pi_1
\end{align}
we get
$$\left((\delta_{a_0}\circ\dots\circ\delta_{a_k})(\Delta_1)\right)(\pi_1 \varphi)=0.$$
Since  $\pi_1$ is surjective it follows that
$$(\delta_{a_0}\circ\dots\circ\delta_{a_k})(\Delta_1)=0,$$
i.e. $\Delta_1\in \Diff_k(A_1)$.
\end{proof}
\newpage

\begin{lm}
Let  $\Delta_i\in \Diff_k(A_i)$, $i=1,2$ be two differential operators. Assume there exists a linear map  $\Delta\colon A\to A$, such that the following diagram commutes
\be\vxymatrix{
A_1\ar[d]^{\Delta_1}&A\ar[d]^\Delta\ar@{->>}[l]_{\pi_1}\ar@{->>}[r]^{\pi_2}&A_2\ar[d]^{\Delta_2}\\
A_1                 &A             \ar@{->>}[l]_{\pi_1}\ar@{->>}[r]^{\pi_2}&A_2
}\ee{d4}
Then such  $\Delta$ is unique and belongs to $\Diff_k(A)$.
\end{lm}
\begin{proof}
We need to prove that $\Delta\in \Diff_k(A) $, i.e. that $$(\delta_{a_0}\circ\dots\circ\delta_{a_k})(\Delta)=0$$
$$(\delta_{a_0}\circ\dots\circ\delta_{a_k})(\Delta)=\sum_{I,\bar{I}} a_I \circ \Delta \circ\ a_{\bar{I}},$$ 
where operator $a_I$ denotes multiplying by elements of algebra, $$I\subset \{0,\dots k\}, \bar{I}=\{0,\dots k\}\setminus I.$$ 
Consider $b\in A$. It is easy to see that $\pi_1(\sum a_I \circ \Delta \circ\ a_{\bar{I}})(b)=0$
\end{proof}
Let us state the following fact using  notation of lemmas  2.2 and 2.3.
\begin{thm}
$$\Diff_{k}(A)=\{(\Delta_1, \Delta_2) \mid \Delta_i\in\Diff_{k}(A_i), p_1\circ\Delta_1\circ\pi_1 = p_2\circ\Delta_2\circ\pi_2 \}.$$ 
\end{thm}
\begin{proof}
It is left to recall how we constructed algebra $A$ given algebras $A_i$. We know that  $p_1(a_1)=p_2(a_2)$, where $a_i\in A_i.$
Hence the following condition holds:
$$p_1\circ\Delta_1\circ\pi_1 = p_2\circ\Delta_2\circ\pi_2,$$ q.e.d.
\end{proof}
\newpage
\section{\Large One dimentional manifolds with singularities}
\subsection{\Large Algebra of functions on a manifold with singularity.}
Let us consider the situation of two curves on a plane which have contact of order $m$. Let us show that the algorithm developed above applies here. 
Let us stick to the case of when the curves are graphs of two functions $h_1(x), h_2(x)$ which have single common point $x=0$. 
Each graph $\Gamma_i$ is a closed subset of the plane hence $\Spec_{\mbbR}C^{\infty}(\Gamma_i) = \Gamma_i$ $[8]$.

Therefore the alebras are $A_i=C^{\infty}(\Gamma_i)$. Moreover by specifying the algebra  $C$ and morphisms  $p_i$ we get a Cartesian square. 
Hence all the algebraical constructions apply.

Now let us show that the algebra of functions  on  $\mathbb{K}_m=\Gamma_1\cup\Gamma_2$ (in the sence of restriction of algebra of functions on  $\mathbb{R}^2$) consides with the algebra of pairs
of functions $f$ and $g$ such that $f(x)-g(x)=o(x^{m})$.
\begin{lm}
$${\mathbb{C}^{\infty}(\mathbb{R}^2)|}_{\mathbb{K}_m} = \{(f,g)\mid f(x)-g(x)=o(x^m)\}.$$
\end{lm}
\begin{proof}

Let us describe $\mathbb{K}_m$ as the subset of plane with coordinate system $(x,y)$. For convenience let us choose coordinates so that
$y(y-h(x))=0$, where  $h(x)$ is a function which has only one zero of order  $(m+1)$ at the point  $x=0$.
Given two functions $f$ and $g$ we can explicitly describe $F(x,y)$: $$F(x,y)= f(x)+ \frac{y(g(x)-f(x))}{h(x)}.$$

Given $F(x,y)$ it is easy to check that $F(x,0)- F(x,h(x))=o(x^m).$ (For instance, by using Hadamard's lemma.)
\end{proof}
Therefore we have a purely algebraical construction for the algebra of smooth functions $C^\infty(\mathbb{K}).$ Similar considerations may be applied 
to other types of singularities.

\newpage
\subsection{\Large Singularity of non-zero order}
\subsubsection{\Large Description of pullback.}
Let us consider two curves on the plane which have a contact of order $m$. According to lemma  3.1 one can apply algorithm described in sections 2.2. and 2.3 to this situation. 
Then
$$A_1=A_2=C^\infty(\mbbR), C=\mbbR[x]/(x^{k+1}),$$ $$p_i(f)=\sum_{n=0}^m \frac{f^{(n)}(0) \epsilon^n}{n!};$$
Algebra of smooth functions on the union of curves may be described as:
$$A=C^\infty(\mathbb{K}_{m})= \{(f(x),g(y))\mid f^{(i)}(0)=g^{(i)}(0), i=0,\dots ,m\}.$$
Algebra $A$ is the Cartesian square as in the second example. Hence results of the previous sections may be applied to it.

\subsubsection{\Large Differential operators}
Let us employ theorem 2.4.
\\
$\Diff_{k}(A)=\{(\Delta_1, \Delta_2)\ |\ $ are such that the following conditions hold$\}:$
\\
$$\Delta_1 =\sum_{i=0}^k  a_i(x){\left(\frac{\partial}{\partial x}\right)}^{i},$$
\\
$$\Delta_2 =\sum_{i=0}^k  b_i(x){\left(\frac{\partial}{\partial y}\right)}^{i}.$$
\\
Conditions on the coefficients follow from the definition, i.e. we have the equality $$p_1\circ\Delta_1\circ\pi_1 =
p_2\circ\Delta_2\circ\pi_2.$$
\\
Since we know the explicit form of $p_i$ the condition transforms into:
$${\Delta_1(f)}^{(i)}(0)={\Delta_2(g)}^{(i)}(0), i=0,\dots ,m.$$
\newpage
For $i\leq m$ we get the following conditions (here  $i$ is the order of the derivative):
\\
$i=0$
\\
$a_s(0)=b_s(0), s\leq m$
\\
$a_s(0)=0=b_s(0), s>m$
\\
\dots
\\
$i=l$
\\
${a_s}^{(l)}(0)={b_s}^{(l)}(0), s\leq m$
\\
${a_s}^{(l)}(0)=0={b_s}^{(l)}(0), s>m+l$
\\
$$\sum_{r=0}^l \binom{l}r a_{m+t-r}^{(l-r)}(0) = 0 = \sum_{r=0}^l \binom{l}r b_{m+t-r}^{(l-r)}(0),\quad\text{для}\quad t=1,\dots ,l,\quad\text{а}\quad l\leq m.$$
Let us assume that  $a_s = 0$ if $s<0$ and $s>k.$
\\
For instance, in case of  $m=1$, considering all the relations on the coeffitients we may write differential operators as:
$$\Delta_x = x^2 \nabla_x -c {\partial_x}^2 + (cx+d)\partial_x + ex+f,$$
$$\Delta_y = y^2 \nabla_y -c {\partial_y}^2 + (cx+d)\partial_y + ey+f,$$
where  $c,d,e,f \in \mathbb{R}$, а $\nabla_x,\nabla_y$ \,--- are some differential operators.

\subsubsection{\Large Symbol of  order-$k$}
Let us consider  $\Delta\in\Diff_k(A)$, it corresponds to a pair of differential operators
$\Delta_1\in\Diff_k(A_1),\Delta_2\in\Diff_k(A_2)$, which satisfy
the conditions of theorem 2.4. Then  $$[\Delta_1]_{k} = a_k(x){\left(\frac{\partial}{\partial
x}\right)}^{k}\quad\text{и}\quad[\Delta_2]_{k} = b_k(y){\left(\frac{\partial}{\partial
y}\right)}^{k}. $$
\\
Using the conditions on the coefficients of differential operators and together with  $smbl_l(\Delta)\cdot smbl_n(\nabla)\in \Smbl_{l+n}$ 
we can define $k$-th symbol via:$$[\Delta]_{k}=([\Delta_1]_{k},[\Delta_2]_{k}),$$
moreover $a_k(0)=0=b_k(0)$, if  $k\geq m$ и $a_k(0)=b_k(0)$, if $k<m$ conditions on the derivatives also hold.

For example, for the contact of order one the conditions would look as follows:
\\
$a_0(0)=b_0(0)$
\\
$a_k(0)=0=b_k(0)$, if $k>0$ 
\\
$a'_k(0)=b'_k(0),$ if $k\leq 2$
\\
$a'_k(0)=0=b'_k(0)$, if $k>2$ 
\newpage
\subsubsection{\Large Poisson bracket}
Consider  $\Delta\in\Diff_l(A)$ and $\nabla\in\Diff_n(A)$, where $\Delta=(\Delta_1,\Delta_2)$ and $\nabla=(\nabla_1,\nabla_2)$, and
$\Delta_i,\nabla_i$ satisfy theorem  2.4. Denote by $a_{i_{\Delta}},b_{i_{\Delta}}$ the coeffitients of  $\Delta$ and by 
$a_{i_{\nabla}},b_{i_{\nabla}}$ the coefficients of $\nabla$. Then  $$\{[\Delta],[\nabla]\}=(\{[\Delta_1],[\nabla_1]\},\{[\Delta_2],[\nabla_2]\}),$$ where
$$\{[\Delta_1],[\nabla_1]\}=(la_{l_{\Delta}}(x)a'_{n_{\nabla}}(x)-na_{n_{\nabla}}(x)a'_{l_{\Delta}}(x)){\left(\frac{\partial}{\partial x}\right)}^{l+n-1},$$
$$\{[\Delta_2],[\nabla_2]\}=(lb_{l_{\Delta}}(y)b'_{n_{\nabla}}(y)-nb_{n_{\nabla}}(y)b'_{l_{\Delta}}(y)){\left(\frac{\partial}{\partial
y}\right)}^{l+n-1}.$$

We get the following conditions on the coeffitients:
\\
$a_i(0)=b_i(0),$ if $i<m$ 
\\
$a'_i(0)=b'_i(0)$, if $i\leq m$
\\
$a_i(0)=0=b_i(0),$ if $i\geq m$ 
\\
$a'_i(0)=0=b'_i(0)$, if  $i>m$  
\\
where $i=l_{\Delta}, n_{\nabla}$
\newpage
\subsection{\Large Singularity of order zero}
\subsubsection{\Large Description of the pullback.}
Let us consider coordinate cross on the plane:
$$\mathbb{K}_{0} = \{(x,y)\subset {\mathbb{R}}^2\mid xy=0\}.$$
Algebra of smooth functions on the cross is given by the formula:
$$ A=C^\infty(\mathbb{K}_{0})={C^{\infty}({\mathbb{R}}^2)|}_{\mathbb{K}_{0}}= \{(f(x),g(y))| f(0)=g(0)\}$$
Applying the algorithm we get in the above notations
$$A_1=A_2=C^\infty(\mbbR), C=\mbbR , p_1(f)=f(0),\,p_2(g)=g(0).$$
\subsubsection{\Large Differential operators}
Theorem 2.4 implies
$$\Diff_{k}(A)=\{(\Delta_1, \Delta_2)\mid a_0(0)=b_0(0),a_i(0)=0=b_i(0),i=1,\dots,k \}.$$
$$\Delta_1 =\sum_{i=0}^k  a_i(x){\left(\frac{\partial}{\partial x}\right)}^{i},$$
$$\Delta_2 =\sum_{i=0}^k  b_i(x){\left(\frac{\partial}{\partial y}\right)}^{i}.$$
\subsubsection{\Large Symbol of order $k$}
Let us use the knowledge we have about contact of non-zero order.
$$[\Delta]_{k}=([\Delta_1],[\Delta_2]),$$ plus
$a_0(0)=b_0(0)$
\\
$a_k(0)=0=b_k(0)$,if $k\geq 1$ 

\subsubsection{\Large  Poisson bracket}
$$\{[\Delta],[\nabla]\}=(\{[\Delta_1],[\nabla_1]\},\{[\Delta_2],[\nabla_2]\}),$$
where $$\{[\Delta_1],[\nabla_1]\}=(la_{l_{\Delta}}(x)a'_{n_{\nabla}}(x)-na_{n_{\nabla}}(x)a'_{l_{\Delta}}(x)){\left(\frac{\partial}{\partial x}\right)}^{l+n-1},$$
$$\{[\Delta_2],[\nabla_2]\}=(lb_{l_{\Delta}}(y)b'_{n_{\nabla}}(y)-nb_{n_{\nabla}}(y)b'_{l_{\Delta}}(y)){\left(\frac{\partial}{\partial
y}\right)}^{l+n-1}.$$
We get the following conditions on the coefficients:
\\
$a_0(0)=b_0(0),$ $a_i(0)=0=b_i(0),$ если $i\geq 1,$ где $i =l_{\Delta}, n_{\nabla}.$

\subsubsection{\Large Spectrum of the algebra of symbols}
Denote by $Smbl(A)=S_0\bigoplus S_1 \bigoplus \dots$ any algebra of symbols.
\begin{lm}For the coordinate cross $|\Smbl(A)|$ is
0-dimentional at the point of singularity and one-dimentional at other points.
\end{lm}
\begin{proof} Let us define homomorphism from the algebra of symbols to $\mbbR.$ 
$$\xymatrix{Smbl(A)=Smbl_0(A)\bigoplus Smbl_1(A)\bigoplus \dots\ar@{.>}[d]^-{H}\\ \mbbR}$$
$H|_{S_0}=h_0,$ we know that  $S_0=A$ $[8].$ Since  $SpecA=\mathbb{K}_{0},$ we known how $H|_{S_0}$ looks. 
Let $h_0$ be homomorphism of evaluation at zero. Let us consider  $(a_k,b_k)\in \Smbl_k,$  
\\
by Hadamards lemma  $$(a_k,b_k)=(a'(0)x+x^2\tilde{a},b'(0)y+y^2\tilde{b})=(a'(0)x,b'(0)y)+(x,y)(x\tilde{a},y\tilde{b}),$$ 
where $$(a'(0)x,b'(0)y)\in\Smbl_k, (x,y)\in\Smbl_0, (x\tilde{a},y\tilde{b})\in\Smbl_k, \tilde{a},\tilde{b}\in C^{\infty}(\mbbR)$$ 
$$H_k(a,b)=H_k(a'(0)x,b'(0)y), k>0.$$
Since $a'(0)=0=b'(0),$  $H_k(a,b)=0.$
$$(H_k(a,b))^2=H_{2k}((a,b)(a,b))=H_{2k}(a^2,b^2)=0,$$
therefore, $H_k(a,b)=0.$
\\
When $h_0$ is evaluation at an arbitrary point let us consider a point $z\in M.$ If $h_o\in |A|$ it may be considered as a homomorphism which maps  $f\in A$ to $f(z).$
\end{proof}
\subsection{\Large Singularity of the first order}
\subsubsection{\Large Description of the pullback}
Let us consider two curves on the plane with contact of order one. By lemma  3.1 algorithm described in sections  2.2. and 2.3 applies. Hence
$$A_1=A_2=C^\infty(\mbbR), C=\mbbR[x]/(x^2),$$ $$p_i(f)=\sum_{n=0}^1 \frac{f^{(n)}(0) \epsilon^n}{n!};$$
Algebra of smooth functions on the union of curves may be described as:
$$A=C^\infty(\mathbb{K}_{1})= \{(f(x),g(y))\mid f^{(i)}(0)=g^{(i)}(0), i=0,1\}.$$
Alebra  $A$ is a Cartesian square and coincides with example 2, hence results of the previous sections apply.

\subsubsection{\Large Differential operators}
By theorem  2.4
$$\Diff_{k}(A)=\{(\Delta_1, \Delta_2)\}:$$
$$\Delta_1 =\sum_{i=0}^k  a_i(x){\left(\frac{\partial}{\partial x}\right)}^{i},$$
$$\Delta_2 =\sum_{i=0}^k  b_i(x){\left(\frac{\partial}{\partial y}\right)}^{i}.$$
Moreover 
\\
$\Diff_0(A)=\{(\Delta_{1_0}, \Delta_{2_0})\mid a_0(0)=b_0(0), a'_0(0)=b'_0(0) \}.$
\\
$\Diff_1(A)=\{(\Delta_{1_1}, \Delta_{2_1})\mid a_0(0)=b_0(0), a'_0(0)=b'_0(0), 
\\
a_1(0)=0=b_1(0), a'_1(0)=b'_1(0) \}.$ 
\\
$\Diff_2(A)=\{(\Delta_{1_2}, \Delta_{2_2})\mid a_0(0)=b_0(0), a'_0(0)=b'_0(0), 
\\
a_1(0)=b_1(0), a'_1(0)=b'_1(0),a_1(0)=b_1(0),a'_1(0)=b'_1(0), 
\\
a_2(0)=0=b_2(0), a'_2(0)+a_1(0)=0=b'_2(0)+a_1(0)\}.$
\\
\dots
\\
$\Diff_l(A)=\{(\Delta_{1_l}, \Delta_{2_l})\mid \dots a_l(0)=0=b_l(0), a'_l(0)=0=b'_l(0), l=3\dots k\}.$

\subsubsection{\Large Symbol of $k$-th order}
Using the above let us define symbols via:
\\
$\Smbl_0 = \{(a, b)\mid a(0)=b(0), a'(0)=0=b'(0)\}.$
\\
$\Smbl_1 = \{(a, b)\mid  a(0)=0=b(0), a'(0)=b'(0)\}.$
\\
$\Smbl_2 = \{(a, b)\mid a(0)=0=b(0), a'(0)=b'(0)\}.$
\\
\dots
\\
$\Smbl_l = \{(a, b)\mid a(0)=0=b(0), a'(0)=0=b'(0), l=3\dots k \}.$

\subsubsection{\Large Poisson bracket}
$$\{[\Delta],[\nabla]\}=(\{[\Delta_1],[\nabla_1]\},\{[\Delta_2],[\nabla_2]\}),$$
where $$\{[\Delta_1],[\nabla_1]\}=(la_{l_{\Delta}}(x)a'_{n_{\nabla}}(x)-na_{n_{\nabla}}(x)a'_{l_{\Delta}}(x)){\left(\frac{\partial}{\partial x}\right)}^{l+n-1},$$
$$\{[\Delta_2],[\nabla_2]\}=(lb_{l_{\Delta}}(y)b'_{n_{\nabla}}(y)-nb_{n_{\nabla}}(y)b'_{l_{\Delta}}(y)){\left(\frac{\partial}{\partial
y}\right)}^{l+n-1}.$$
Conditions on the coefficients are:
\\
$a_0(0)=b_0(0)$  
\\
$a'_i(0)=b'_i(0)$, if $i\leq 2$
\\
$a_i(0)=0=b_i(0),$ if $i\geq 1$ 
\\
$a'_i(0)=0=b'_i(0)$, if $i>2$  
\\
where $i=l_{\Delta}, n_{\nabla}$

\subsubsection{\Large Spectrum of algebra of symbols}
Algebra of symbols is  $Smbl(A)=S_0\bigoplus S_1 \bigoplus \dots$
\begin{lm}For the curves with contact of order one  $|\Smbl(A)|$ is 0-dimentional at singular point and 1-dimentional at the others.
\end{lm}
\begin{proof}
 Let us define the homomorphism from alebra of symbols to  $\mbbR$ by
$$\xymatrix{\Smbl(A)=\Smbl_0(A)\bigoplus \Smbl_1(A)\bigoplus \dots\ar@{.>}[d]^-{H}\\R}.$$
\\
Let $H_0=h_0$, i.e. evaluation at zero. Let us consider $(a_k,b_k)\in \Smbl_k,$  
by Hadamards lemma  $$(a_k,b_k)=(a'(0)x+x^2\tilde{a}),b'(0)y+y^2\tilde{b})=(a'(0)x,b'(0)y)+(x,y)(x\tilde{a},y\tilde{b}),$$ 
where $$(a'(0)x,b'(0)y)\in\Smbl_k, (x,y)\in\Smbl_0, (x\tilde{a},y\tilde{b})\in\Smbl_k.$$
\\
Since $a'(0)=b'(0)=c,$ we have
$$H_1(a,b)=cH_1(a,b)+H_1(x^2a''(0),y^2b''(0)),$$
$$H_2(a,b)=cH_2(a,b)+H_2(x^2a''(0),y^2b''(0)),$$
$$H_k(a,b)=H_k(x^2a''(0),y^2b''(0)), k>2.$$
Let  $a=x,b=y$, then  $H_3(x^3,y^3)=0=(H_1(x,y))^3,$ hence $H_1(a,b)=0,$ similarly one may prove that  $H_2(a,b)=0.$
\\
Therefore we found that  $H_k(a,b)=0,$ i.e. spectrum at the point of singularity is a set of one point.
\end{proof}
\newpage
\section{\Large Conclusion}
To develop Hamiltonian formalizm on a manifold one needs to define the Poisson bracket.
Well known results of classical theory belonging to Poisson Hamilton, Ostrogradskii and Liouville were obtained for 
a bracket defined in canonical coordinates in the phase space.
Algebraic definition of Poisson bracket is presented in this paper. Algebra $A$ with the induced Poisson bracket is called Poisson manifold. Hamiltonian is the element of algebra $A$ .
Hamiltonian system is the triple $(A,\{,\},H)$.

Correspondence between a hamiltonian and a vector field is given by the formula
$X_H(f)=\{H,f\}$. It has the following physical meaning: a point of the phase space moves along the field $X_H$.

Algebraic approach allows to develop a unified formalism for the classical and quantum mechanics.


\end{document}